\documentclass[12pt]{article}
\usepackage{amsmath,amssymb,amsthm,fullpage}

\title{Quantum symmetry groups of noncommutative spheres}

\author{Joseph C. V\'arilly\footnotemark \\
        Departamento de Matem\'atica,
        Universidad de Costa Rica\\
        2060 San Jos\'e, Costa Rica}

\date{19 June 2001}

\renewcommand{\a}{\alpha}         %% group action
\renewcommand{\b}{\beta}          %% group action
\newcommand{\Dl}{\Delta}          %% coproduct
\newcommand{\dl}{\delta}          %% Dirac delta
\newcommand{\eps}{\varepsilon}    %% counit
\newcommand{\ga}{\gamma}          %% group action
\newcommand{\La}{\Lambda}         %% exterior power
\newcommand{\la}{\lambda}         %% constant
\newcommand{\sg}{\sigma}          %% cocycle
\renewcommand{\th}{\theta}        %% torus parameter

\newcommand{\A}{\mathcal{A}}      %% algebra
\renewcommand{\H}{\mathcal{H}}    %% Hilbert space

\newcommand{\R}{\mathbb{R}}       %% real numbers
\newcommand{\Sf}{\mathbb{S}}      %% sphere
\newcommand{\T}{\mathbb{T}}       %% circle as a group
\newcommand{\Z}{\mathbb{Z}}       %% integers

\newcommand{\g}{\mathfrak{g}}     %% Lie algebra
\newcommand{\hl}{\mathfrak{h}}    %% abelian Lie algebra

\DeclareMathOperator{\ad}{ad}     %% infml adjoint action
\DeclareMathOperator{\Aut}{Aut}   %% automorphism group
\DeclareMathOperator{\ch}{ch}     %% Chern character
\DeclareMathOperator{\id}{id}     %% identity map

\newcommand{\hookto}{\hookrightarrow} %% inclusion arrow
\newcommand{\op}{\oplus}          %% direct sum
\newcommand{\ox}{\otimes}         %% tensor product
\newcommand{\x}{\times}           %% Moyal or direct product
\renewcommand{\.}{\cdot}          %% anonymous variable
\renewcommand{\:}{\colon}         %% colon in  f: A -> B

\newcommand{\Coo}{C^\infty}       %% smooth functions
\newcommand{\Dslash}{{D\mkern-11.5mu/\,}} %% Dirac operator
\newcommand{\thalf}{\tfrac{1}{2}} %% small fraction  1/2
\DeclareMathOperator{\tsum}{{\textstyle\sum}} %% small sum in display
\newcommand{\hatox}{\mathrel{\widehat\otimes}} %% proj tensor product
\newcommand{\row}[3]{{#1}_{#2},\dots,{#1}_{#3}}   %% list: a_1,...,a_n
\newcommand{\sepword}[1]{\quad\text{#1}\quad} %% well-spaced words
\newcommand{\twobytwoeven}[2]{\begin{pmatrix}#1 & 0 \\
                            0 & #2\end{pmatrix}} %% 2 x 2 diag matrix
\newcommand{\avg}[1]{\langle#1\rangle} %% average or VEV
\theoremstyle{plain}
\newtheorem*{thm}{Theorem}

\theoremstyle{definition}
\newtheorem*{ack}{Acknowledgments}

\theoremstyle{remark}
\newtheorem{exmp}{Example}

\begin{document}

\maketitle

\renewcommand{\thefootnote}{\fnsymbol{footnote}}
\addtocounter{footnote}{1}
\footnotetext{Regular Associate of the Abdus Salam ICTP, Trieste;
email: \texttt{varilly@cariari.ucr.ac.cr}.}

\begin{abstract}
We show that the noncommutative spheres of Connes and Landi are
quantum homogeneous spaces for certain compact quantum groups.
We give a general construction of homogeneous spaces which support
noncommutative spin geometries.
\end{abstract}

\section{Introduction}
\label{sec:intro}

Noncommutative geometry \cite{Book} has established itself as a theory
which goes beyond the realm of differentiable manifolds and deals in a
unified fashion with many singular geometric spaces, too. A
fundamental feature of NCG is that it fully incorporates all compact,
boundaryless spin manifolds under the heading of ``noncommutative spin
geometries'': see \cite{ConnesGrav} and Chapter 11 of \cite{Polaris}.

Outstanding examples of singular geometric spaces are the
noncommutative tori \cite{ConnesTorus,ConnesToral,RieffelRot}, orbit
spaces of discrete group actions, and leaf spaces of foliations.
Recently, a new class of examples has appeared, the ``noncommutative
spheres'' of Connes and Landi~\cite{ConnesLa}, from a purely
cohomological construction.

The Moyal-like nature of the twisted products introduced in
\cite{ConnesLa} suggests that the underlying noncommutative
spaces of these spin geometries may be obtained, as $C^*$-algebras, by
the general deformation construction of Rieffel~\cite{RieffelDefQ}.
The question arises as to whether these are in fact noncommutative
homogeneous spaces, that is, subalgebras of invariants of certain Hopf
algebras which may be regarded as ``quantized symmetry groups''. This
question is more delicate than it might seem, because it must be
answered at the $C^*$-algebra level: these ``symmetry groups'' must be
found in the category of ``compact quantum groups'' in the sense of
Woronowicz~\cite{WoronowiczComQGp} or perhaps in the wider category of
``locally compact quantum groups'' \cite{KustermansVLocComp}. As it
happens, the compact noncommutative spaces which we discuss below have
compact (quantum) symmetry groups, so we shall restrict ourselves here
to Woronowicz' version.

In Sections~\ref{sec:quant-iv-sf} and~\ref{sec:def-homog} we review
the construction of noncommutative spheres and Rieffel's
$C^*$-deformation theory. Section~\ref{sec:comp-qgp-def} treats
compact quantum groups built by such deformations. In
Section~\ref{sec:ncsf-homog}, we explain how both constructions mesh
to yield the desired quantum homogeneous spaces. In the final section,
we briefly discuss noncommutative spin geometries on these homogeneous
spaces.

\section{Quantized 4-spheres}
\label{sec:quant-iv-sf}

The construction of noncommutative spin geometries by Connes and Landi
proceeds in two stages. First, the data $(\A,\H,D,C,\chi)$ of an even
real spectral triple~\cite{Book,Polaris} are sought as possible
solutions to a system of equations for the Chern character in cyclic
homology:
\begin{subequations}
\label{eq:ch-eqns}
\begin{align}
\ch_k(p) \equiv \avg{(p - \thalf)\,dp^{2k}} 
&= 0  \sepword{for} k = 0,1,\dots,m-1,
\label{eq:ch-vanish}  \\
\pi_D(\ch_m(p)) & = \chi,
\label{eq:ch-volume}
\end{align}
\end{subequations}
where $p = p^2 = p^*$ is an orthogonal projector in a matrix algebra
$M_r(\A)$, $\avg{\.}$ denotes the conditional expectation (or partial
trace) onto $\A$, $\chi$ is the grading operator on the $\Z_2$-graded
Hilbert space $\H$, and
$\pi_D(a_0\,da_1\dots da_n) := a_0\,[D,a_1]\dots [D,a_n]$ represents
elements of the universal graded differential algebra over $\A$ as
operators on~$\H$.

These equations impose restrictions, first of all, on the algebra $\A$
itself. In dimension two, i.e., when $m = 1$ and $r = 2$, only
commutative solutions are found; in fact, Connes showed by an
elementary argument \cite{ConnesSurvey} ---see also
\cite[Sect.~11.A]{Polaris} and \cite{Paschke}--- that
\eqref{eq:ch-vanish} alone forces $\A$ to be a commutative algebra
whose Gelfand spectrum is a closed subset of the 2-sphere $\Sf^2$.
This equation also makes $\ch_1(p)$ a Hochschild $2$-cycle, whose
associated volume form is the standard volume form on the sphere, so
the Gelfand spectrum must be the whole $\Sf^2$, and thus
$\A \simeq \Coo(\Sf^2)$ on the basis of \eqref{eq:ch-eqns} alone!

Even in commutative cases such as this, where $D$ may be taken as the
Dirac operator given by some metric and spin structure on the spectrum
of~$\A$, the final condition \eqref{eq:ch-volume} does not determine
the metric, but only its volume form; thus the cohomological
conditions \eqref{eq:ch-eqns} allow for volume-preserving variations
of the metric, as befits a theory which aspires to incorporate
gravity.

In dimension four, with $m = 2$ and $r = 4$, there is also a
commutative solution given in~\cite{ConnesSurvey}, namely the smooth
function algebra $\Coo(\Sf^4)$. Later, Connes and Landi
\cite{ConnesLa} found a family of noncommutative solutions,
parametrized by a complex number of modulus one $\la = e^{2\pi i\th}$:
these are the algebras $\Coo(\Sf_\th^4)$ (together with their
corresponding Dirac operators), which may be called ``smooth function
algebras for noncommutative $4$-spheres $\Sf_\th^4$'', in the standard
parlance of quantum group theorists. Their representations are
uniformly bounded and in each case a $C^*$-norm is quickly found,
allowing to complete them to ``continuous function algebras'', denoted
$C(\Sf_\th^4)$.

This procedure extends directly to higher dimensions, yielding
noncommutative spheres in any even dimension greater than~$2$ from the
corresponding ``instanton algebras'' (so called because the finite
projective modules $p\A^r$ may be regarded as vector bundles
over~$\A$). Starting from the odd Chern character in cyclic homology,
one can also search for odd-dimensional noncommutative spaces with
this method (in the odd case, $\H$ is ungraded and $\chi$
in~\eqref{eq:ch-volume} is replaced by~$1$).

A striking feature of this construction is that these noncommutative
manifolds are parametrized by numbers of modulus one, in contrast
to the \textit{real} numbers $q \neq \pm 1$ which label the well-known
$2$-spheres $\Sf_{qc}^2$ of Podle\'s \cite{Podles}, which were
originally constructed as homogeneous spaces of the compact quantum
groups $SU_q(2)$. By combining features of both constructions,
D\c{a}browski, Landi and Masuda \cite{DabrowskiLM} built a family of
quantized $4$-spheres $\Sf_q^4$; on computing the Chern characters of
the instantons, they found that \eqref{eq:ch-vanish} is violated,
inasmuch as $\ch_1(p) = (1 - q^2)$ times a nonvanishing term.

In any case, it is clear that the Connes--Landi spheres $\Sf_\th^4$
lie outside the realm of $q$-spheres of the Podle\'s type. Indeed,
several other variants on the $\Sf_q^4$ spheres have since appeared
\cite{BonechiCT,BrzezinskiG,SitarzSphere}, which, however, do not
incorporate the $\Sf_\th^4$ family \cite{DabrowskiL}. Of particular
note is the construction by Hong and Szyma\'nski~\cite{HongS} of a
large family of quantized $n$-spheres $\Sf_q^n$, for $n \geq 2$ and
$q > 0$, by deforming $C(\Sf^n)$ to Cuntz--Krieger $C^*$-algebras
based on certain directed graphs; but again, the $\Sf_\th^4$ family is
not included. Therefore, it behoves us to ask whether that family may
be realized as ``quantum homogeneous spaces''.

\section{Deformations of homogeneous spaces}
\label{sec:def-homog}

The second stage of the Connes--Landi construction is the provision of
spin geometries on the spheres $\Sf_\th^4$. This is accomplished by a
deformation of the commutative spectral triple
$(\Coo(\Sf^4),\H,\Dslash)$, where $\Dslash$ denotes a Dirac operator
on the Hilbert space $\H$ of square-integrable spinors over $\Sf^4$.
In the deformation, $\Dslash$ is kept fixed, so that all spectral
data, including the classical dimension (four!) of the geometry are
unchanged: only the algebra and its representation on $\H$ are
modified.

One declares a kind of Moyal product on
$\Coo(\Sf^4)$ by the following recipe: first, note that there is an 
isometric action of the $2$-torus $\T^2$ on $\Sf^4$, allowing us to 
decompose any smooth function on $\Sf^4$ as a series
$f = \sum_r f_r$ indexed by $r \in \Z^2$, where $f_r$ lies in the 
$r^{\mathrm{th}}$ spectral subspace:
$$
(e^{2\pi i\phi_1}, e^{2\pi i\phi_2}) \. f_r
 = e^{2\pi i(r_1\phi_1 + r_2\phi_2)} \, f_r.
$$
The series converges rapidly in the Fr\'echet topology of
$\Coo(\Sf^4)$. By introducing the following star-product of
homogeneous elements:
\begin{equation}
f_r \x g_s := e^{2\pi i\th r_1s_2}\, f_r g_s,
\label{eq:CoLa-starprod}
\end{equation}
Connes and Landi constructed a representation of $C(\Sf_\th^4)$ on the
spinor space~$\H$ (having bounded commutators with $\Dslash$); in
essence, the representation is explicit only on the smooth subalgebra,
which is just the vector space $\Coo(\Sf^4)$ with the commutative
product replaced by the star-product~\eqref{eq:CoLa-starprod}.

More generally, if $M$ is a compact Riemannian manifold admitting a
Lie group of isometries of rank $l \geq 2$, so that $M$ carries an
isometric action of the torus $\T^l$, one can decompose $\Coo(M)$ into
spectral subspaces indexed by $\Z^l$. The Moyal product of two 
homogeneous functions $f_r$ and $g_s$ is then given by
\begin{equation}
f_r \x g_s := \rho(r,s)\, f_r g_s,
\label{eq:Moyal-rs}
\end{equation}
where $\rho\: \Z^l \x \Z^l \to \T$ is a $2$-cocycle on the additive
group $\Z^l$. The cocycle relation
$$
\rho(r, s + t) \rho(s,t) = \rho(r,s) \rho(r + s, t)
$$
guarantees associativity of the new product. For instance
\cite{Polaris,Atlas}, one may take 
$$
\rho(r,s) := \exp\bigl\{-2\pi i \tsum_{j<k} r_j \th_{jk} s_k \bigr\},
$$
where $\th = [\th_{jk}]$ is a real $l \x l$ matrix. Complex
conjugation of functions remains an involution for the new product
provided that the matrix $\th$ is \textit{skewsymmetric}.

The relation \eqref{eq:Moyal-rs} is easily recognized as the
product rule for the twisted group $C^*$-algebra $C^*(\Z^l,\rho)$. We 
may replace $\rho$ by its skewsymmetrized version
\begin{equation}
\sg(r,s) := \exp\bigl\{-\pi i \tsum_{j,k=1}^l r_j \th_{jk} s_k\bigr\},
\label{eq:tori-cocyc}
\end{equation}
because $\rho$ and $\sg$ are cohomologous~\cite{RieffelPMod}, and we
obtain $C^*(\Z^l,\sg) = C(\T_\th^l)$, which is precisely the
$C^*$-algebra of the noncommutative $l$-torus with parameter
matrix~$\th$.

\vspace{6pt}

The relation \eqref{eq:Moyal-rs} is clearly, then, a discretized
version of the usual Moyal product, due to the periodicity of the
$\T^l$-action. Recall that the standard Moyal product on the phase
space $\R^{2m}$ may be expressed either by the familiar series in
powers of~$\hbar$ whose first nontrivial term gives the Poisson
bracket, or alternatively in the integral form~\cite{Phobos}:
$$
(f \x_J g)(x) := (2\pi\hbar)^{-n}
 \iint f(x + s) g(x + t) \,e^{is\.Jt/\hbar} \,ds\,dt,
$$
where $J$ is the skewsymmetric matrix giving the standard symplectic
structure on $\R^{2m}$ (and the dot is the usual scalar product
on~$\R^{2m}$). This may be interpreted as an oscillatory integral for
suitable classes of functions and distributions on $\R^{2m}$, and
yields the familiar series as an \textit{asymptotic} expansion in
powers of~$\hbar$ \cite{Nereid,Voros}. It is, therefore, a better
starting point than that series for a $C^*$-algebraic theory of
deformations. Indeed, this was the form of the Moyal product used by
Rieffel in his general deformation theory~\cite{RieffelDefQ}. He
found, in fact, an improvement over the board by rewriting it as
$$
(f \x_J g)(x) := \iint f(x + Js) g(x + t) \,e^{2\pi is\.t} \,ds\,dt.
$$
He then generalized this to
\begin{equation}
a \x_J b
 := \iint_{V\x V} \a_{Js}(a) \a_t(b) \,e^{2\pi is\.t} \,ds\,dt,
\label{eq:Marc-prod}
\end{equation}
where $a,b$ belong to a $C^*$-algebra $A$, $\a\: V \to \Aut(A)$ is a
(strongly continuous) action of a vector group $V \simeq \R^l$ on~$A$,
and $J$ is a skewsymmetric real $l \x l$ matrix. The oscillatory
integral \eqref{eq:Marc-prod} makes sense, \textit{a priori}, only for
elements $a,b$ of the smooth subalgebra $A^\infty$ of~$A$ (under the
action~$\a$), which is a Fr\'echet pre-$C^*$-algebra.

This problem of good definition is overcome \cite{RieffelDefQ} by
introducing a suitable $C^*$-norm on $A^\infty$ for which the $\x_J$
product is continuous, and then completing it in this norm to obtain
the deformed $C^*$-algebra $A_J$. The construction is functorial in
that morphisms of~$A$ restrict to $A^\infty$ and then extend uniquely
to morphisms of~$A_J$. In more detail: if $(A,\a(V))$ and $(B,\b(V))$
are two $C^*$-algebras carrying actions of~$V$, and if
$\phi\: A \to B$ is a $*$-homomorphism intertwining the actions $\a$
and~$\b$, then $\phi(A^\infty) \subseteq B^\infty$ and the restriction
of $\phi$ to $A^\infty$ extends uniquely to a $*$-homomorphism
$\phi_J \: A_J \to B_J$. Moreover, if the original map $\phi$ is
injective, then $\phi_J$ is injective, too; and $\phi_J$ is surjective
whenever $\phi$ is surjective.

In particular, when $B = A$ and $\b = \a$, each $\a_x$ intertwines
$\a$ with itself since $V$ is an abelian group, and this gives an
action $\a_J \: V \to \Aut(A_J)$ whose restriction to $A^\infty$
coincides with $\a$. Then $(A_J,\a_J)$ can be deformed in turn, using
a new skewsymmetric matrix $K$, say; and the result turns out to be
isomorphic to $A_{J+K}$. By taking $K = -J$, we see that the change
$A \mapsto A_J$ is reversible. It is therefore unsurprising, but still
a deep and important result, that the smooth subalgebra remains
unchanged during this mutation: $(A_J)^\infty = A^\infty$ as vector
spaces, although they have different multiplications
\cite[Thm.~7.1]{RieffelDefQ}.

The case of particular interest to us occurs when the action $\a$ of 
$V$ is periodic, so that $\a_x = \id_A$ for $x \in L$, a cocompact 
lattice in~$V$; in which case, $\a$ is effectively an action of the 
compact abelian group $H = V/L$. Then $A^\infty$ decomposes into 
spectral subspaces labelled by elements of $L$ (or characters of~$H$) 
and one can check \cite[Prop.~2.21]{RieffelDefQ} that if 
$\a_s(a_p) = e^{2\pi ip\.s} a_p$ and $\a_t(b_q) = e^{2\pi iq\.t} b_q$
with $p,q \in L$, then 
$$
a_p \x_J b_q = e^{-2\pi ip\.Jq} a_p b_q.
$$
On comparing this with \eqref{eq:Moyal-rs} (with the cocycle $\rho$ 
replaced there by~$\sg$), we see that it suffices to take 
$A := C(\T^l)$ and $J := \thalf\th$ in order to obtain any 
noncommutative torus $C(\T_\th^l) \simeq A_J$ by this algorithm.

In fine, the isospectral deformation procedure of \cite{ConnesLa},
based on the star-product~\eqref{eq:CoLa-starprod}, is, as far as the
algebra is concerned, a special case of Rieffel's $C^*$-deformation
theory. (This same point is made by Sitarz in a recent
announcement~\cite{SitarzDef}.)

Moreover, if the isometric action of $\T^l$ on~$M$ alluded to above is
free on some orbit, so that $M$ contains an embedded $l$-torus, then
the restriction map $\pi\: C(M) \to C(\T^l)$ induces a
surjective $*$-homomorphism $\pi_{2J}\: C(M)_{2J} \to C(\T_\th^l)$, so
that the noncommutative torus appears as a quotient of the deformed
$C(M)$. In particular, if $\th$ is irrational, then the noncommutative
sphere $C(\Sf_\th^4)$ is not a type~I $C^*$-algebra.

\section{Compact quantum groups from deformations}
\label{sec:comp-qgp-def}

It stands to reason, then, that this $C^*$-deformation process should
yield compact quantum groups when applied to the $C^*$-algebra $C(G)$
of continuous functions on a compact Lie group. This proves to be the
case, by a further construction of Rieffel. There are two issues to
address here: first, which vector group actions on $C(G)$ are
admissible and useful, and second, how to deal with the coproduct,
counit and antipode which define the Hopf algebra structure of $C(G)$
(or rather, of its dense subalgebra of representative functions).

The solution to the second problem could not be simpler: the coalgebra
structure and antipode can be left completely untouched, and only the
algebra structure need be deformed! The matter is not quite trivial,
as one must ensure that the coproduct is still an algebra homomorphism
for the new product. This possibility was pointed out by
Dubois-Violette \cite{DuboisVQGp}, who noticed that Woronowicz' matrix
corepresentations for $C(SU_q(N))$ and similar bialgebras could be
seen as different star-products on the same coalgebra.

Now suppose that $H$ is a closed connected abelian subgroup of~$G$
(usually we may take $H$ to be a maximal torus, but it is not really
necessary that it be maximal); following Rieffel~\cite{RieffelComQGp},
we consider the action of $H \x H$ on~$G$ given by
$(h,k) \. x := h x k^{-1}$, and the corresponding action on $C(G)$:
\begin{equation}
[(h,k) \. f](x) := f(h^{-1} x k).
\label{eq:HGH-action}
\end{equation}
We may regard this as a periodic action of the Lie algebra
$\hl \op \hl$, with the following notation. Choose and fix a basis for 
the vector space $\hl \simeq \R^l$, so that the exponential mapping 
from $\hl$ onto $H$ may be expressed as a homomorphism
$e\: \R^l \to H$ whose kernel is the integer lattice $\Z^l$; by taking
$\la := e(1,1,\dots,1)$ we may write $\la^s := e(s)$ for $s \in \R^l$
with a multiindex notation; the action of $V := \hl \op \hl$ on $C(G)$
is then written as
\begin{equation}
[\a(s,t) f](x) := f(\la^{-s} x \la^t).
\label{eq:hCGh-action}
\end{equation}
(In the sequel, we shall refer to this as an action of $\hl \op \hl$
or of $H \x H$, interchangeably.) If $J$ is now any skewsymmetric
matrix in $M_{2l}(\R)$, then \eqref{eq:Marc-prod} now defines a Moyal
product on $\Coo(G)$, and the procedure of Sect.~\ref{sec:def-homog}
extends this to a $C^*$-algebra $C(G)_J$, which yields a quantization
of~$G$ as a noncommutative space.

The remaining difficulty is that an arbitrary choice of~$J$ will not
mesh well with the coalgebra structure of $\Coo(G)$, so we shall not
always get a quantization of the \textit{group} structure of~$G$. For
that, one may follow the approach of Drinfeld by first equipping $G$
with a compatible Poisson bracket (i.e., the product map
$G \x G \to G$ must be a Poisson map). By a well-known procedure
\cite{ChariP} this can be done at the infinitesimal level by equipping
its Lie algebra $\g$ with a cocycle $\phi\: \g \to \La^2\g$ whose dual
defines a Lie bracket on~$\g^*$. For instance, one may take
$\phi(X) = \ad_X(r)$, where $r \in \g \ox \g$ is a skewsymmetric
solution of the classical Yang--Baxter equation
$[r_{12}, r_{13}] + [r_{12}, r_{23}] + [r_{13}, r_{23}] = 0$. If 
$r = \sum_k X_k \ox Y_k$, then since
$[r_{12}, r_{13}] = \sum_{jk} [X_j,X_k] \ox Y_j \ox Y_k$ and similarly 
for the other terms, this equation is satisfied when
$r \in \hl \ox \hl$ for an \textit{abelian} Lie subalgebra $\hl$
of~$\g$. Although there are other solutions (see \cite{LevendorskiiS}
for an exhaustive treatment of Poisson Lie group structures on
\textit{simple} compact Lie groups and the several algebraic
quantizations of the Hopf algebra of representative functions), we
shall focus on the case $r \in \hl \ox \hl$. On using our previous
identification of $\hl$ with~$\R^l$, we can write $r$ as a
skewsymmetric $l \x l$ matrix~$Q$. The corresponding Poisson structure
on~$G$ is given by the bivector field $W$, where
$W_x := \la_x(r) - \rho_x(r)$ is the difference of the left and right
translates of $r$ from $\La^2\g$ to $\La^2 T_xG$.
Therefore~\cite{RieffelNonComp}, at the infinitesimal level we should
take
$$
J := \twobytwoeven{Q}{-Q}
$$
as the $2l \x 2l$ matrix of deformation parameters for the action of 
$\hl \op \hl$.

We can now write the twisted product on $\Coo(G)$ as 
\begin{equation}
(f \x_J g)(x) := \int_{\hl^4} f(\la^{-Qs} x \la^{-Qt})
 g(\la^{-u} x \la^v) \,e^{2\pi i(s\.u + t\.v)} \,ds \,dt \,du \,dv.
\label{eq:Moyal-gpprod}
\end{equation}

The coproduct $\Dl$, the counit $\eps$ and the antipode $S$, which are
defined on the Hopf algebra of representative functions of~$G$ by
\begin{equation}
\Dl f(x,y) := f(xy),  \qquad  \eps(f) := f(1),  \qquad 
Sf(x) := f(x^{-1}),
\label{eq:coalg-relns}
\end{equation}
whereby $\Dl$ and $\eps$ are algebra homomorphisms and $S$ is an
antiisomorphism, obviously extend to algebra maps of $C(G)$ with the
same properties. It is shown in \cite{RieffelComQGp,Wang} that they
also satisfy the same algebraic relations for the twisted product. The
formulas \eqref{eq:coalg-relns} make sense for $f \in \Coo(G)$ or even
$f \in C(G)$, although the usual requirement
$\Dl(\Coo(G)) \subseteq \Coo(G) \ox \Coo(G)$ holds only if the
algebraic tensor product is replaced by the completed tensor product,
which we denote by $\Coo(G) \hatox \Coo(G)$ and identify with
$\Coo(G \x G)$. We can make a formal check of these homomorphism
properties for smooth functions:
\begin{align*}
&(\Dl f \x_J \Dl g)(x,y)
\\
&= \int_{\hl^8} f(\la^{-Qs} x \la^{-Qt-Qs'} y \la^{-Qt'})
     g(\la^{-u} x \la^{v-u'} y \la^{v'})
     \,e^{2\pi i(s\.u + t\.v + s'\.u' + t'\.v')} \,ds \dots dv'
\\
&= \int_{\hl^8} f(\la^{-Qs} x \la^{-Qt''} y \la^{-Qt'})
     g(\la^{-u} x \la^{-u''} y \la^{v'})
     \,e^{2\pi i(s\.u + t''\.v + s'\.u'' + t'\.v')} \,ds \dots dv'
\\
&= \int_{\hl^6} f(\la^{-Qs} x \la^{-Qt''} y \la^{-Qt'})
     g(\la^{-u} x \la^{-u''} y \la^{v'})
     \,e^{2\pi i(s\.u + t'\.v')} \,\dl(t'') \,\dl(u'') \,ds \dots dv'
\\
&= \int_{\hl^4} f(\la^{-Qs} xy \la^{-Qt'}) g(\la^{-u} xy \la^{v'})
     \,e^{2\pi i(s\.u + t'\.v')} \,ds\,dt'\,du\,dv'
\\
&= (f \x_J g)(xy) = \Dl(f \x_J g)(x,y).
\end{align*}
Similarly,
\begin{align*}
(f \x_J g)(1)
&= \int_{\hl^4} f(\la^{-Q(s+t)}) g(\la^{v-u})
     \,e^{2\pi i(s\.u + t\.v)} \,ds \,dt \,du \,dv
\\
&= \int_{\hl^4} f(\la^{-Q(s')}) g(\la^{v'})
     \,e^{2\pi i(s'\.u + t\.v')} \,ds' \,dt \,du \,dv'
\\
&= \int_{\hl^2} f(\la^{-Q(s')}) g(\la^{v'}) \,\dl(s') \,\dl(v')
     \,ds' \,dv' = f(1) \, g(1),
\end{align*}
so $\eps(f \x_J g) = \eps(f) \eps(g)$. Next, if $Q$ is invertible, then
\begin{align*}
(Sf \x_J Sg)(x)
&= \int_{\hl^4} f(\la^{Qt} x^{-1} \la^{Qs}) g(\la^{-v} x^{-1} \la^u)
     \,e^{2\pi i(s\.u + t\.v)} \,ds \,dt \,du \,dv
\\
&= (\det Q)^{-2} \int_{\hl^4} f(\la^{-t'} x^{-1} \la^{s'})
     g(\la^{-v} x^{-1} \la^{-u})
     \,e^{-2\pi i(Q^{-1}t'\.v + Q^{-1}s'\.u)} \,ds' \,dt' \,du \,dv
\\
&= \int_{\hl^4} f(\la^{-t'} x^{-1} \la^{s'}) 
     g(\la^{-Qv'} x^{-1} \la^{-Qu'}) 
     \,e^{2\pi i(t'\.v' + s'\.u')} \,ds' \,dt' \,du' \,dv'
\\
&= (g \x_J f)(x^{-1}) = S(g \x_J f)(x),
\end{align*}
using the skewsymmetry of $Q$ in the third step; on the other hand, if
$Q = 0$, then $f \x_J g = fg$ and the calculation reduces to
$(Sf \x_J Sg)(x) = f(x^{-1}) g(x^{-1}) = S(g \x_J f)(x)$; since we may
integrate separately over the nullspace of~$Q$ and its orthogonal 
complement, we conclude that $Sf \x_J Sg = S(g \x_J f)$ in all cases.

The defining property of the antipode may also be checked in this
manner: indeed, if $m(f \ox g) := f \x_J g$ for $f, g \in \Coo(G)$,
similar formal calculations quickly establish that
\begin{equation}
m(\id \ox S)(\Dl f) = \eps(f)\,1 = m(S \ox \id)(\Dl f)
\label{eq:antp-defn}
\end{equation}
whenever $f \in \Coo(G)$. However, it should be pointed out that the
previous calculations in fact involve oscillatory integrals of
functions of $s,t,u,v \in \R^l$ which have neither compact support nor
fast decrease; but with some additional careful analysis, it is shown
in \cite{RieffelDefQ} that they remain valid for smooth functions
which have all derivatives bounded on~$\R^l$, as is always the case
when $f,g \in \Coo(G)$.

In summary, the product \eqref{eq:Moyal-gpprod} on $\Coo(G)$ is fully
compatible with its original coalgebra structure and antipode. The
functoriality of the $A_J$ construction then lifts $\Dl$ and $S$ as
algebra (anti)homomorphisms to the $C^*$-level. (Some bookkeeping is
necessary because the source and target algebras carry different
actions of $\hl \op \hl$ in each case.) With $A = C(G)$ and
$J = Q \op (-Q)$ as before, we then obtain a continuous
$*$-homomorphism $\Dl_J \: A_J \to A_J \ox A_J$ (with the minimal
$C^*$-tensor product) and a continuous $*$-antihomomorphism
$S_J \: A_J \to A_J$. The counit $\eps$ on $\Coo(G)$ also extends to a
character of~$A_J$.

Note, however, that the twisted product on $\Coo(G)$ generally does
not extend to a continuous linear map from $A_J \ox A_J$ to $A_J$.
(For one thing, $m$ is not an algebra homomorphism unless $G$ is
abelian.) Thus, the relation \eqref{eq:antp-defn} is not helpful at
the $C^*$-level. This is an old problem, and for unital 
$C^*$-algebras there is a well-known solution, described in the 
fundamental paper of Woronowicz~\cite{WoronowiczComQGp}. Given a 
\textit{unital} $C^*$-algebra $A$ and a unital $*$-homomorphism 
$\Dl\: A \to A \ox A$ which is coassociative, define linear maps $W$,
$W'$ on the algebraic tensor product of $A$ with itself by
$$
W(a \ox b) := (\Dl a)(1\ox b) \sepword{and}
W'(a \ox b) := (a\ox 1)(\Dl b).
$$
(These are the Kac--Takesaki or ``fundamental unitary'' operators.)
Woronowicz' postulate is that the maps $W$, $W'$ have dense range.
Then $(A,\Dl)$ is called a \textit{compact quantum group}. The counit
and antipode are automatically defined on a dense $*$-subalgebra, and
$A$ has a unique state (the ``Haar state'') which is both left and
right invariant \cite{WoronowiczComQGp}. For $A = C(G)$, these maps
are
$$
W(f \ox g)(x,y) := f(xy)g(y),  \qquad  W'(f \ox g) := f(x)g(xy),
$$
which have dense range in $C(G \x G)$. After deformation, these become
$$
W(f \ox g)  := (\Dl f) \x_J (1 \ox g),  \qquad
W'(f \ox g) := (f \ox 1) \x_J (\Dl g),
$$
for $f,g \in \Coo(G)$, and these extend to invertible maps on 
$\Coo(G \x G)$. Concretely, for $h \in \Coo(G \x G)$,
\begin{align*}
Wh(x,y)
&= \int_{\hl^4} h(x \la^{-Qs} y \la^{-Qt}, \la^{-u} y \la^v)
     \,e^{2\pi i(s\.u + t\.v)} \,ds\,dt\,du\,dv,
\\
W^{-1}h(x,y)
&= \int_{\hl^4} h(x \la^{Qt} y^{-1} \la^{-Qs}, \la^{-u} y \la^v)
     \,e^{2\pi i(s\.u + t\.v)} \,ds\,dt\,du\,dv,
\end{align*}
as may be verified directly. It follows that $W$ and likewise $W'$ 
have dense range in $C(G \x G)$.

\section{Noncommutative spheres as homogeneous spaces}
\label{sec:ncsf-homog}

The standard $4$-sphere is a homogeneous space of the $5$-dimensional
rotation group: $\Sf^4 \approx SO(5)/SO(4)$. Note that $SO(5)$ is a
compact simple Lie group of rank two. More generally, we may consider
homogeneous spaces of the form $M = G/K$, where $G$ is a compact Lie
group (which need not be semisimple) and $K$ is a closed subgroup. Let
$H$ be a closed abelian subgroup of~$K$; then we can deform both
$C(G)$ and $C(K)$ by the \textit{same} action~\eqref{eq:hCGh-action}
of $H \x H$. Note in passing that a maximal torus in $SO(2l)$ is
carried onto a maximal torus of $SO(2l+1)$ by the standard inclusion
$SO(2l) \subset SO(2l+1)$, so that even-dimensional spheres
$\Sf^{2l} = SO(2l+1)/SO(2l)$ fall under this heading.

The left action of $G$ on $G/K$ yields a $*$-homomorphism
$\rho\: C(G/K) \to C(G) \ox C(G/K)$ by $\rho f(x,yK) := f(xyK)$.
Restricted to smooth functions, this can be viewed as a left coaction
of $\Coo(G)$ on $\Coo(G/K)$. Let $C(G)^K$ denote the subalgebra of
$C(G)$ consisting of right-invariant functions under the action of~$K$,
so $f \in C(G)^K$ if $f(xw) = f(x)$ whenever $w \in K$, $x \in G$; and
let $\Coo(G)^K := C(G)^K \cap \Coo(G)$. There is an obvious
$*$-isomorphism $\zeta\: C(G)^K \to C(G/K)$ given by
$\zeta f(xK) := f(x)$, and $\zeta(\Coo(G)^K) = \Coo(G/K)$. The
coproduct $\Dl$ of $\Coo(G)$ maps $\Coo(G)^K$ into
$\Coo(G) \hatox \Coo(G)^K$, the space of smooth functions $h$ on
$G \x G$ for which $h(x,yw) \equiv h(x,y)$ when $w \in K$. Moreover,
if $f \in \Coo(G)^K$, then
$$
[\rho\zeta f](x,yK) = \zeta f(xyK) = f(xy) = \Dl f(x,y)
 = [(\id \ox \zeta)\Dl f](x,yK),
$$
so $\zeta$ intertwines the coactions $\rho$ and $\Dl$. In short, the
algebra $\Coo(G/K)$, together with its isomorphism onto $\Coo(G)^K$,
is an embedded homogeneous space in the Hopf algebra $\Coo(G)$.

Now we come to the main point. Since $H \subseteq K$, the left-right 
action \eqref{eq:HGH-action} of $H \x H$ on both $G$ and $K$ induces a 
left action of $H$ on $G/K$, since the right action of~$H$ is absorbed 
in the right $K$-cosets. If we deform $C(G)$ and $C(K)$ via the 
$H \x H$ action along the direction $J = Q \op (-Q)$, the 
corresponding effect on $C(G/K)$ should be a deformation under an
$H$-action along the direction~$Q$. And so it proves. 

To see that, we first notice that for $f,g \in \Coo(G)^K$,
\eqref{eq:Marc-prod} yields
\begin{align*}
(f \x_J g)(x) 
&= \int_{\hl^4} f(\la^{-Qs} x \la^{-Qt}) g(\la^{-u} x \la^v)
    \,e^{2\pi i(s\.u + t\.v)} \,ds\,dt \,du\,dv
\\
&= \int_{\hl^4} f(\la^{-Qs} x) g(\la^{-u} x)
    \,e^{2\pi i(s\.u + t\.v)} \,ds\,dt \,du\,dv
\\
&= \int_{\hl^2} f(\la^{-Qs} x) g(\la^{-u} x) \,e^{2\pi is\.u}\,ds\,du,
\end{align*}
or
\begin{equation}
f \x_J g = \int_{\hl^2} \ga_{Qs}(f) \ga_u(g) \,e^{2\pi is\.u}\,ds\,du
\label{eq:xJ-xQ}
\end{equation}
where $(\ga_t f)(x) := f(\la^t x)$ for $f \in C(G)^K$. The action of
$\hl$ on $C(G/K)$ may be defined as $(\b_t h)(xK) := h(\la^t xK)$, so
that $\zeta$ intertwines the actions $\b$ and $\ga$ of~$H$. Then
\eqref{eq:xJ-xQ} becomes simply
$$
\zeta f \x_Q \zeta g = \zeta(f \x_J g)
  \sepword{for all}  f,g \in \Coo(G/K).
$$

Finally, we can lift this isomorphism to the $C^*$-level, using the
functoriality of $C^*$-deformations. First, since
$\zeta\: C(G/K) \to C(G)^K$ is a $*$-isomorphism intertwining $\b$ and
$\ga$, its restriction to $\Coo(G/K)$ extends to a $*$-isomorphism of
$C(G/K)_Q$ to $C(G)^K_Q$, where the latter comes from the action of
$\ga$ on $\Coo(G)^K$. Of course, $\ga$ can be regarded as an action of
$H \x H$ where the second factor acts trivially; since elements of
$C(G)^K$ are right-invariant under $H$, $\ga$ is just the restriction
of the action $\a$ to~$C(G)^K$. This means that the inclusion
$C(G)^K \hookto C(G)$ is equivariant for the actions $\ga$ and $\a$,
and so its restriction to $\Coo(G)^K$ extends to a $*$-homomorphism
from $C(G)^K_Q$ to $C(G)_J$; by Proposition~5.8 of~\cite{RieffelDefQ},
this is still injective. In summary,
$$
C(G/K) \simeq C(G)^K \hookto C(G)  \sepword{leads to} 
C(G/K)_Q \simeq C(G)^K_Q \hookto C(G)_J.
$$

If the subgroup $H$ is not a maximal torus in either $K$ or~$G$, the
space of smooth elements for the action of $H \x H$ will be strictly
larger than $\Coo(G)$ (for instance, if the action is trivial, all
continuous functions are smooth in this sense); however, as clarified
in Sect.~1 of~\cite{RieffelComQGp}, we may continue to use $\Coo(G)$
instead, because it will be dense in the Fr\'echet topology of the
space of all smooth elements, and therefore will remain dense in the
deformed $C^*$-algebra $C(G)_J$. The same applies, \textit{mutatis 
mutandis}, to $\Coo(G/K)$ and $C(G/K)_Q$.

\vspace{6pt}

We have thus proved the following result.

\begin{thm}
The deformed $C^*$-algebra $C(G/K)_Q$ is an embedded 
homogeneous space for the compact quantum group $C(G)_J$.
\qed
\end{thm}

\begin{exmp}
The even-dimensional noncommutative spheres $\Sf_\th^{2l}$ of Connes 
and Landi come directly from this framework, for $l \geq 2$.
Just take $G = SO(2l+1)$, $K = SO(2l)$ and let $H \simeq \T^l$ be a 
maximal torus for~$K$; then let $Q = \thalf\th$, where $\th$ is a 
skewsymmetric $l \x l$ matrix.

The odd-dimensional spheres $\Sf^{2l+1} = SO(2l+2)/SO(2l+1)$ have
somewhat different deformations, since the $l$-dimensional maximal
torus of $SO(2l+1)$ is not maximal in $SO(2l+2)$, so the twisted
product reduces to the ordinary commutative product along some 
directions.
\end{exmp}

\begin{exmp}
Our construction yields several new examples of homogeneous spaces.
For instance, if $T$ is a maximal torus of~$G$, the flag manifold
$G/T$ may be deformed in any direction $Q = -Q^t$ in $M_l(\R)$
provided $l = \dim T \geq 2$. In particular, it yields a family of
$6$-dimensional quantized manifolds $C(SU(3)/\T^2)_Q$. It would be
of interest to classify these up to isomorphism or Morita equivalence.

At the algebraic level, there are other deformations of flag 
manifolds~\cite{LevendorskiiS} which go beyond those considered here, 
in that more general solutions of the classical Yang--Baxter equation
are used for the deformation directions. These could yield further 
examples of quantum homogeneous spaces.
\end{exmp}

\section{Homogeneous noncommutative spin geometries}
\label{sec:homog-NCgeom}

These new homogeneous spaces give rise to spectral triples, by the
isospectral deformation procedure of~\cite{ConnesLa}. We may start
from the manifold $G/K$ with, say, the normalized $G$-invariant
metric. Suppose that $G/K$ also has a homogeneous spin structure (if
not, a homogeneous spin$^\mathrm{c}$ structure will do). Let $D$ be
the corresponding Dirac operator, let $\row{X}{1}{l}$ be the chosen
basis of $\hl$, and let $p_j$ be the selfadjoint operator representing
$X_j$ on the spinor space $\H$, for $j = 1,\dots,l$. Since the action
of~$\hl$ integrates to a representation of~$H$ on spinors, the
operators $p_j$ have integer or half-odd-integer spectra, and for each
$r \in \Z^l$, there is a unitary operator
$\sg(p,r) := \exp\bigl\{-2\pi i \tsum_{j,k} p_j Q_{jk} r_k\bigr\}$,
using the notation of~\eqref{eq:tori-cocyc}; its inverse is
$\sg(r,p)$. These operators commute with each other and also with $D$,
although not with the representation of $\Coo(G/K)$ on~$\H$. Any
bounded operator $T$ in the common smooth domain of the
transformations $T \mapsto \sg(p,r) T \sg(r,p)$ has a decomposition
$T = \sum_{r\in\Z^l} T_r$, where $\sg(p,r)\,T_s = T_s\,\sg(p+s,r)$
for $r,s \in \Z^l$; define
$$
L(T) := \sum_{r\in\Z^l} T_r \,\sg(p,r).
$$
The cocycle property of~$\sg$ immediately gives
$L(f)L(g) = L(f \x_Q g)$, so that $L$ yields a representation of
$(\Coo(G/K), \x_Q)$ on~$\H$, while
$[D, L(f)] = \sum_r [D,f_r]\, \sg(p,r) = L([D,f])$ is a bounded
operator for all $f \in \Coo(G/K)$. The charge conjugation operator
$C$ on spinors~\cite[Chap.~9]{Polaris} commutes with all $\sg(p,r)$
and therefore $C p_j C^{-1} = -p_j$ for each~$j$. It follows that
$R(T) := C L(T)^* C^{-1}$ is given by
$$
R(T) = \sum_{r\in\Z^l} \sg(r,p)\, C T_r^* C^{-1}
     = \sum_{r\in\Z^l} C T_r^* C^{-1} \,\sg(r,p).
$$
Since $C f^* C^{-1} = f$ for $f$ in the commutative algebra
$\Coo(G/K)$, this reduces to $R(f) = \sum_{r\in\Z^l} f_r \,\sg(r,p)$,
and therefore $R(f)R(g) = R(f \x_{-Q} g)$. (Our use of the
skewsymmetrized cocycle $\sg$ obviates the need to twist the
conjugation as in~\cite{ConnesLa}.) It is easy to see ---compare
\cite{Phobos}--- that $R$ gives an antirepresentation of
$(\Coo(G/K), \x_Q)$ on~$\H$, which commutes with~$L$ because
\begin{align*}
L(f) R(g)
&= \sum_{r,s} f_r \,\sg(p, r) \,g_s \,\sg(s, p)
 = \sum_{r,s} f_r g_s \,\sg(p+s, r) \,\sg(s, p)
\\
&= \sum_{r,s} g_s f_r \,\sg(s, p+r) \,\sg(p, r)
 = \sum_{r,s} g_s \,\sg(s, p) \,f_r \,\sg(p, r) = R(g) L(f).
\end{align*}
This verifies the reality property of the spin geometry. It is readily
checked that
$$
[[D, L(f)], R(g)]
 = \sum_{r,s\in\Z^l} \sg(p,r)\, [[D, f_r], g_s] \,\sg(s,p) = 0,
$$
so the first-order property of the spin geometry holds, too.

Such spin geometries $(L(\Coo(G/K)),\H,D,C,\chi)$ has maximal
symmetry; we may indeed refer to the quantum group $C(G)_J$ as its
``noncommutative symmetry group''. They provide examples of spectral
triples with noncommutative symmetries as discussed, for instance,
in~\cite{PaschkeS}. However, only the invariance of~$D$ under the
abelian subgroup $H$ is actually used, so we are free to build other
spin geometries by deforming the commutative ones obtained from any
$H$-invariant metric on $G/K$.

More elaborate examples of deformed geometries can also be built,
starting from commutative spin geometries wherein the spin connection 
is replaced by a Clifford superconnection (as in~\cite{Zappafrank}, 
for instance), provided the latter is also $H$-invariant.

\vspace{6pt}

Finally, we consider whether the noncommutative homogeneous spaces
constructed here may play the same role as noncommutative tori in
quantum field theory. Recall that Seiberg and
Witten~\cite{SeibergWGeom} and Konechny and Schwarz~\cite{KonechnyS}
have extensively explored noncommutative gauge theories based on tori.
In general, the divergent ultraviolet behaviour for field theories
based on noncommutative tori~\cite{Atlas} is no better than in the
commutative case. This divergence holds also for field theories
obtained by second-quantizing the spin geometries constructed here.

Without going into the detailed analysis, the matter may be summed up
as follows. The action of a $G$-invariant Dirac operator over $G/K$
decomposes into matrix actions on finite-dimensional subspaces of
smooth spinors, for which explicit formulas are
available~\cite{Baer,Slebarski}. The sign operator $F := D|D|^{-1}$
preserves these subspaces, which are permuted by the representation of
the algebra $(\Coo(G/K),\x_Q)$. For any unitary $u$ in this algebra,
we can decompose the operator $[F,u]$ as in~\cite{Atlas} or
\cite[Sect.~13.A]{Polaris} and estimate its Schatten class, which
measures the degree of ultraviolet divergence of the theory. The norms
$\|[F,u]\|_p$ turn out to be independent of the cocycle $\sg$ defining
the product, provided $\sg(r, r + s) = \sg(r,s)$; in view
of~\eqref{eq:tori-cocyc}, this is immediate from the skewsymmetry of
the parameter matrix~$Q$. Therefore, the overall UV behaviour remains
the same as in the commutative case when $Q = 0$: our deformations 
never soften the ultraviolet divergence.

The ubiquity of the Moyal product in noncommutative field theory is
already familiar. While the present work cannot pretend to explain its
pervasiveness, we have at any rate shown that noncommutative
geometries with a high degree of symmetry are easy to deform along (at
least two) commuting directions, leading always to Moyal products with
a few parameters; thus the emphasis on noncommutative tori is by no
means misplaced. Whether this is in the nature of things remains to be
seen.

\begin{ack}
We thank A. Connes, L. D\c{a}browski, H. Figueroa,
J.~M. Gracia-Bond\'{\i}a, P.~M. Hajac and M. Paschke for helpful
discussions on several matters. Support from the Vicerrector\'{\i}a de
Investigaci\'on of the University of Costa Rica and the Abdus Salam
ICTP, Trieste, is gratefully acknowledged.
\end{ack}

\end{document}